\documentclass[12pt]{article}

\newlength{\stefan}
\usepackage{amsmath, amsthm, amsfonts, makeidx}
\usepackage{multind}
\DeclareMathSymbol{\subsetneq}{\mathord}{AMSb}{"26}

\newtheorem{lemma}{Lemma}[section]

\newtheorem{theorem}[lemma]{Theorem}

\theoremstyle{definition}

\newtheorem{remark}[lemma]{Remark}
\newtheorem{conjecture}[lemma]{Conjecture}

\newcommand{\lp}{\longrightarrow}

\newcommand{\F}{\mathbb{F}}

\newcommand{\E}{\mathcal{E}}

\newcommand{\kk}{\mathbb{K}}

\newcommand{\GA}{\textup{GA}}

\newcommand{\sym}{\operatorname{Sym}}
\newcommand{\alt}{\operatorname{Alt}}
\newcommand{\GL}{\textup{GL}}

\title{A problem on polynomial maps over finite fields}
\author{Stefan Maubach\\
\small Dept. of Math.\\
\small Radboud University Nijmegen\\
\small The Netherlands\\
\small  s.maubach@math.ru.nl\\
\small http://www.math.ru.nl/$\sim$maubach}

\begin{document}

\maketitle

Recently I noticed that an old conjecture of mine got quite some interest of people working outside of my field,
especially people who work in discrete mathematics, number theory, code theory, combinatorics, even cryptographers.
Due to that interest, I thought it might be a nice idea to write a short instructive
article giving a very down-to-earth explanation for someone completely unfamiliar with it. Also, 
I give the background of the problem and indicate its importance.
This all with the idea that the problem enters the world of discrete mathematics, and can be approached from completely different angles than the ones I, and my colleagues,
are familiar with.

Well, first let me state the question without explaining the notations.\\

\noindent
{\bf Problem:} Do there exist odd polynomial automorphisms over the finite fields $\F_4,\F_8,\F_{16},\F_{32},\ldots$?\\

\noindent
Below I will explain all terms used in this problem, in rather elaborate terms. You might want to skip some parts, depending on your knowledge. \\

\section{Polynomial automorphisms}

For simplicity, I will explain everything in dimension 2. 
Let $\kk$ be a field.
If $F_1,F_2\in \kk[X,Y]$ then we write $F:=(F_1(X,Y),F_2(X,Y))$ is a polynomial endomorphism.
One can see such a map in several ways:\\
\begin{enumerate}
\item As an element of $(\kk[X,Y])^2$, this is how I explained it above.
\item As a map $\kk^2\lp\kk^2$, induced by polynomials. The map sending $(a,b)\in \kk^2$ to $(F_1(a,b),F_2(a,b))\in \kk^2$.
\item As a map $\kk[X,Y]\lp \kk[X,Y]$, sending $p(X,Y)$ to $p(F_1(X,Y),F_2(X,Y))$.
\end{enumerate}

Let's stick to point 1 for the moment.
We say that $I:=(X,Y)$ is the identity map.
If $F,G$ are polynomial maps, we can make their composition $F\circ G:=(F_1(G_1,G_2),F_2(G_1,G_2))$. If $F\circ G=I$, then we say that $G$ is the inverse of $F$.
(Automatically, $G\circ F=I$.) We say that $F$ is a {\em polynomial automorphism}.\\

Note that polynomial automorphisms include all invertible linear maps (or ``nonsingular matrices'' if you like that more).
Just to give you one example of a ``non-trivial'' polynomial automorphism: $(X+Y^2,Y)$ has as inverse $(X-Y^2,Y)$.

\begin{remark} Obviously, if $F$ is a polynomial automorphism, then $F$ induces a bijection $\kk^2\lp \kk^2$. However, the converse is only true if $\kk$ is infinite!!
Since we are talking exactly about the finite field case, we do NOT have this situation here! For example, $(X^p,Y^p)\in (\F_p[X,Y])^2$ is not a polynomial automorphism,
but it does induce a bijection $\F_p^2\lp \F_p^2$.
\end{remark}

\section{What does ``odd'' mean?}

From now on, $\kk=\F_q$, a field with $q$ elements.
If you have a polynomial automorphism $F$, then it induces a map $\E(F): \F_q^2\lp \F_q^2$.
It's not just a map, it is a bijection of $\F_q^2$. Since these are $q^2$ elements, one can see $\E(F)$ as an element of
the symmetric group with $q^2$ elements. In general, in dimension $n$:
\[ \E(F)\in \sym(q^n). \]
Obviously, you can get any endomorphism of $\F_q^2$ induced by an element of $(\F_q[X,Y])^2$ (even in dimension 1 !) and hence any element of $\sym(q^n)$
can be given by an element of $(\F_q[X,Y])^2$. But, what about which bijections are induced by polynomial automorphisms?
This question is partially answered in \cite{Mau01}, and the answer is quite surprising:

\begin{theorem} \label{T1} Let $n\geq 2$.
\begin{enumerate}
\item
If $q$ is odd or $q=2$ then any element of $\sym(q^n)$ can be obtained as the image of a polynomial  automorphism.
\item
If $q=2^m$ where $m\geq 2$ then any element of the subgroup $\alt(q^n)$ can be obtained as the image of a polynomial  automorphism.
\end{enumerate}
\end{theorem}

The group $\alt(q^n)$ is the alternating subgroup of elements which are even. Obviously, the above theorem is screaming the question:
do there exist elements of $\sym(q^n)\backslash \alt(q^n)$ which can be obtained as the image of a polynomial automorphism over $\F_4, \F_8$ etc.?
So let us repeat the problem:\\

\noindent
{\bf Problem:} Do there exist odd polynomial automorphisms over the finite fields $\F_4,\F_8,\F_{16},\F_{32},\ldots$?\\

\section{Why is it important?}

If one can find an odd polynomial automorphism, then it has many consequences!! Such an automorphism will single-handedly kill quite a few conjectures.
But before I can explain why, I have to throw some theory at you (which is easy enough to read before going to bed, or
while having a headache, or both). \\

Let us denote the polynomial automorphism group in dimension $n$ over $\kk$ by $\GA_n(\kk)$.
One of the big problems in my field (affine algebraic geometry) is trying to understand the group $\GA_n(\kk)$.
There are quite a few problems that we can only address when we have understood the automorphism sufficiently enough.

At the moment, we don't even know if we have a generating set for the polynomial automorphism group in dimension 3 and higher!!

We have a reasonable description in dimension 2, which then has as a side effect that many results are only known in dimension 2 and
not beyond!!\\
\ \\
{\bf \large Tame automorphisms}\\

First, let us make some polynomial automorphisms. As we already noticed, all invertible linear maps are polynomial automorphisms. So:
\begin{equation} \GL_n(\kk)\subseteq \GA_n(\kk) . \end{equation}
Then we have a class of very simple polynomial automorphisms. The {\em elementary} ones. They look like this:
\[ (X_1+f(X_2,\ldots,X_n),X_2,\ldots,X_n) \]
where $f\in \kk[X_2,\ldots,X_n]$ arbitrarily. Its inverse is the elementary map where you replace $f$ by $-f$.
If you compose a few of them, you can make any {\em triangular} map
\[ (aX+f(Y,Z),~ bY+g(Z), ~cZ+h) \]
where $f\in \kk[Y,Z], g\in \kk[Z], h\in \kk, a,b,c\in \kk$. (This is only an example in dimension 3, but it suffices to understand the dimension 3 case.)
 Now
\[ \begin{array}{c}
(X,Y,cZ+h)\circ(X, bY+g(Z), Z)\circ (aX+f(Y,Z), Y, Z)= \\
 (aX+f(Y,Z), bY+g(Z), cZ+h) \end{array}\]
and thus any of such triangular maps is invertible. The set of these maps has many names: two are Jonqui\'ere group and Triangular group.
It is denoted by
\begin{equation}  \textup{BA}_n(\kk) \subseteq \GA_n(\kk) \end{equation}
since it resembles the Borel subgroup of the linear maps (the triangular linear maps).
(In fact, $\textup{BA}_n(\kk) \cap \GL_n(\kk)=B_n(\kk)$, the Borel group.)
Now you can make the group generated by $\textup{BA}_n(\kk)$ and $\GL_n(\kk)$, which is called the {\em tame automorphism group:}
\begin{equation} \textup{TA}_n(\kk):=<\textup{BA}_n(\kk),\GL_n(\kk)>. \end{equation}
We call automorphisms tame if they are elements of $\textup{TA}_n(\kk)$ and wild otherwise.

In dimension 2 there is the famous Jung-van der Kulk-theorem, which states among others (for any field $\kk$):
\[ \textup{TA}_2(\kk)=\GA_2(\kk) .\]
This is the reason why we can do so much more in dimension 2: we have a set of generators of the automorphism group, and even quite a nice description of it.\\

Now in dimension three or higher there were some candidate examples of which no one knew if they are tame or not. The most famous is {\em Nagata's map:}
\[ (X-2Y\Delta -Z\Delta^2,Y+Z\Delta, Z) \]
where $\Delta=XZ+Y^2$. (The inverse is $(X+2Y\Delta -Z\Delta^2,Y-Z\Delta, Z)$, by the way.)\\
\ \\
{\bf \large The biggest breakthrough in the last twenty years} in our field was the proof in 2004 of Shestakov-Umirbaev \cite{SU04, SU04a}, which showed
that Nagata's map (and  others) were not tame. That was and is an amazing result, which surprised our community. For these papers, they got the  AMS 2007 Moore Prize  (a  ``best paper award'' given every three years \cite{AMS}). It's quite a nontrivial paper, and has lots of technicalities. For now, it only works in characteristic zero.

Now, in \cite{Mau01} the above theorem \ref{T1} was actually the following:

\begin{theorem} Let $n\geq 2$.
\begin{enumerate}
\item
If $q$ is odd or $q=2$ then $\E(\textup{TA}_2(\F_q))=\sym(q^n)$.
\item
If $q=2^m$ where $m\geq 2$ then $\E(\textup{TA}_2(\F_q))=\alt(q^n)$.
\end{enumerate}
\end{theorem}

In other words, if you find an ``odd'' polynomial automorphism over for example $\F_4$, then it cannot be tame ! 
Checking that a concrete polynomial automorphism is odd, can be easily done by computer or by hand. 
So, that would be a very, very easy proof of an automorphism being non-tame, giving you a great result in a one-page paper.

As you can guess, all examples I and some others have tried so far are even. But, I do not know what to believe! There are nowadays several ways to make automorphisms,
and it is very unclear if not one of these could be odd.

On the other hand, one could try to prove that some classes of polynomial automorphism over $\F_4,\F_8,\ldots$ are even. These are nice, concrete questions which you can
work wit. 

If you are ambitous, you could try to do this for {\em all} polynomial automorphisms, solving the problem. (If you don't believe there are odd polynomial automorphisms.) 
But it is hard to use the input that you have a 
polynomial automorphism in this particular problem! My guess is that one should consider a larger class of polynomial endomorphisms, and prove something about them. 
This approach could lead to something, I think. If there are no odd polynomial automorphisms, that is\ldots

\section{More consequences of an odd automorphism}

There are several conjectures on generating sets of automorphism groups. Some of them would be killed by an odd automorphism. 

\begin{conjecture} $\GA_3(\kk)$ is generated by the automorphisms fixing one variable. 
\end{conjecture}

If $\kk=\F_4,\F_8,\ldots$ then you can show that any automorphism fixing one variable is even. For, let $(F_1(X,Y,Z), F_2(X,Y,Z), Z)$ be such 
an automorphism, then for any value $Z=a\in \kk$, $(F_1(X,Y,a), F_2(X,Y,a))$ must be an automorphism in dimension 2. Here we have the Jung-van der Kulk theorem, so 
this map is even. So, on each surface $\kk^2\times \{a\}$ the map is even. So the total map is even too. 

\begin{conjecture} $\GA_n(\kk)$ is generated by all linearizable automorphisms. 
\end{conjecture}

An automorphism $F$ is linearizable if there exists another automorphism $\varphi$ such that $\varphi^{-1}F\varphi$ is linear. 
Since all linear maps over $\F_4,\F_8,\ldots$ are even, this conjecture would imply that each map is even, as being even or odd is invariant under conjugation.

\end{document}